\documentclass[12pt]{article}
\title{Symplectic Resolutions for Quotient Singularities}
\author{Baohua FU}
\usepackage{amsmath,amssymb,amsthm}
\chardef\bslash=`\\
\newtheorem{Thm}{Theorem}[section]
\newtheorem{Cor}[Thm]{Corollary}
\newtheorem{Lem}[Thm]{Lemma}
\newtheorem{Prop}[Thm]{Proposition}
\newtheorem{Def}{Definition}
\newtheorem{Rque}{Remark}

\def\cit{{\mathbb C}}

\def\qit{{\mathbb Q}}
\def\pit{{\mathbb P}}
\def\zit{{\mathbb Z}}
\def\0{{\mathcal O}}

\begin{document}
\maketitle
\begin{abstract}
Let $X$ be a smooth irreducible complex variety and $G$ a finite subgroup of $Aut(X)$. There is a natural action of $G$ on $T^*X$ preserving the 
canonical symplectic form. We show that if $T^*X/G$ admits a symplectic resolution  $\pi: Z \rightarrow T^*X/G$, then $X/G$ is smooth and $Z$ contains an 
open set isomorphic to $T^*(X/G)$. In the case of $X = \pit^n$ and $G \subset SL(n+1,\cit)$, we give an McKay correspondence type description for
the Hodge numbers of $Z$.
\end{abstract}
\section{Introduction}

Recall that a regular 2-form on a smooth complex algebraic variety is {\em symplectic} if it is closed and non-degenerate at every point.
A {\em resolution} is a projective morphism from a smooth variety onto an algebraic variety $X$ which is an isomorphism outside of
the singular locus of $X$. 
In the pioneering paper \cite{Be}, A.Beauville has proposed and initiated the following notion of symplectic singularities. 
\begin{Def} 
A normal complex algebraic variety $X$ is said to have {\em symplectic singularities} if there 
exists a regular 
symplectic 2-form $\omega$ on $X_{reg}$ such that for any resolution of singularities $\pi: \widetilde{X} \rightarrow X$, the 2-form
$\pi^* \omega$ defined a priori on $\pi^{-1}(X_{reg})$ can be extended to a regular 2-form on $\widetilde{X}$. If furthermore 
the 2-form $\pi^* \omega$ extends to a  symplectic 2-form on the whole of $\widetilde{X}$ for some resolution of $X$, then we say that $X$ admits
a {\em symplectic resolution}.
\end{Def}

 As shown in \cite{Fu} (see also \cite{Ka1} and \cite{Ve}), a resolution of symplectic singularities is 
symplectic if and only if it is crepant. An important class of examples of symplectic singularities is the normalization of  the closure of a nilpotent
 orbit in a semi-simple Lie algebra. For these singularities, the problem of the 
existence of symplectic resolutions has been completely solved in \cite{Fu}.

{\bf Convention:} throughout this note, all varieties are supposed to be quasi-projective, which guarantees the existence of the quotient by a
 finite group. 

Another class of examples of symplectic singularities comes from the following:
\begin{Prop}[\cite{Be}, Prop. 2.4]
Let $V$ be a variety with symplectic singularities, G a finite group of automorphisms of $V$, preserving a symplectic 2-form on 
$V_{reg}$. Then the variety $V/G$ has symplectic singularities.
\end{Prop}

A very particular case is the quotient of a complex vector space $\cit^{2n}$ by a finite group $G$ of symplectic automorphisms. The problem of the 
existence of
symplectic resolutions for such quotient singularities has been studied by D.Kaledin in \cite{Ka1} and  M.Verbitsky in \cite{Ve}. In particular, they 
have shown that the existence of a symplectic resolution for $\cit^{2n}/G$ implies that $G$ is generated by symplectic reflections, i.e. by the elements
$g \in G$ such that the fixed point set by $g$ in $\cit^{2n}$ is of co-dimension 2.   However the complete classification problem
of $G$ such that $\cit^{2n}/G$ admits a symplectic resolution is far from being solved. 

The importance of this problem  comes from the generalized McKay correspondence.  
If a sympletic resolution $\pi: Z \rightarrow \cit^{2n}/G$ exists, then the
McKay correspondence relates the geometry of $Z$ to some  properties of $G$. An example is the correspondence between a ``natural'' basis of
$H_*^{BM}(Z, \qit)$  and the conjugacy classes of $G$ (see \cite{Ka2} and \cite{Ba}), here $H_*^{BM}(Z,\qit)$ is the Borel-Moore homology groups.  

In this note, we study the existence of symplectic resolutions for general quotient symplectic singularities. We prove in section 2 the following:
\begin{Thm}
Let $V$ be a smooth irreducible symplectic variety and $G$ a finite subgroup of $Aut(V)$ which preserves the symplectic form on $V$. 
Suppose that $V/G$ admits  a symplectic resolution.  Then the closed subvariety $\cup_{g \neq 1} Fix(g)$ is either empty or  of 
pure codimension 2 in $V$.
\end{Thm}
As an easy corollary, an isolated quotient symplectic singularity admits a symplectic resolution if and only if it is of dimension 2. 
In section 3, we restrict to the particular case of $V = T^*X$, where $X$ is a smooth irreducible  complex algebraic variety, and $G$ a finite subgroup of 
$Aut(X)$ which acts naturally on $T^*X$ preserving the canonical symplectic structure on $V$. Then using some ideas of Kaledin \cite{Ka1}, we prove
\begin{Thm}
If the quotient $T^*X/G$ admits a symplectic resolution $\pi: Z \rightarrow  T^*X/G$, then $X/G$ is smooth and there exists a Zariski open set $U$ of $Z$ 
which is isomorphic to the total space of the cotangent bundle $T^*(X/G)$ of $X/G$. 
\end{Thm}

The idea of the proof is to study the natural $\cit^*$-action on $T^*X/G$, which lifts on $Z$. The symplectic form $\Omega$ on $Z$ satisfies 
$\lambda^* \Omega = \lambda \Omega$, which makes the arguments in \cite{Fu} (see also \cite{Ka1} and \cite{Nak}) possible.

Applying this to the case of $X = \cit^n$, we get
\begin{Cor}
If $\pi: Z \rightarrow \cit^n \oplus \cit^n/G$ is a symplectic resolution, then $Z$ is simply connected and rational.
\end{Cor}

In section 4, we study McKay correspondence for symplectic resolutions of $T^*X/G$. The key tool is the orbit E-functions (see \cite{Ba}), which
is shown to be equal to the string E-function in {\em loc.cit.}. In the particular case of $X = \pit^n$, we prove the following:
\begin{Thm}
Let $X= \pit^n$ be the projective space and $G$ a finite group in $SL(n+1,\cit)$. Suppose that $\pi: Z \rightarrow T^*X/G $ is a symplectic
resolution. Then the Hodge numbers of $Z$ can be calculated by the following 
$$\sum_{p,q} (-1)^{p+q} h^{p,q}(Z) u^p v^q = (uv)^n \sum_{\{g\}} \sum_i \cfrac{(uv)^{k_i(g)} - 1}{uv - 1},  $$ where
$k_i(g)$ are multiplicities of distinct eigenvalues of $g \in SL(n+1,\cit)$ and $\{g\}$ runs over conjugacy classes in $G$ considered as a subgroup in 
$Aut(\pit^n)$.
\end{Thm}
\begin{Cor}
The Euler number of $Z$ is given by $e(Z) = (n+1)  c(G) $, where $c(G)$ is the number of conjugacy classes of elements in $G$ considered as 
a subgroup in $Aut(\pit^n)$.
\end{Cor}

In the last section, we prove G\"ottesche's formula 
for Hilbert schemes of the cotangent bundle of a curve, which  has also been proved by H.Nakajima \cite{Nak} using another approach.

{\em Acknowledgments:} I want to thank D.Kaledin, C.Walter and J.Wierzba for some helpful discussions. 
 I am especially grateful to A.Beauville and A.Hirschowitz for many interesting discussions and suggestions.
Without their help, this work could never have been done.

\section{General case}

Recall that a variety is {\em $\qit$-factorial} if any Weil divisor has a multiple that is a Cartier divisor. A smooth variety is $\qit$-factorial. 
\begin{Thm} \label{thm}
Let $V$ be a smooth irreducible symplectic variety and $G$ a finite group of symplectic automorphisms.. 
Suppose that $V/G$ admits  a symplectic resolution.  Then the closed subvariety $F = \cup_{g \neq 1} Fix(g)$ is either empty or  of 
pure  codimension 2 in $V$.
\end{Thm}
\begin{proof}
Being a quotient of a $\qit$-factorial normal variety by a finite group,  $V/G$ is again $\qit$-factorial and normal. This gives that any component
$W$ of the exceptional set of a resolution $\pi: Z \rightarrow V/G$ is of pure codimension one. On the other hand, by a result of D. Kaledin 
(\cite{Ka3} proposition 1.2, see also \cite{Nam} corollary 1.5 for the same assertion without assuming projectivity) any symplectic
 resolution is semismall, i.e. $2 codim(W) \geq codim(\pi(W))$, thus $2 = 2 codim(W) \geq codim(\pi(W))$. Suppose that $V/G$ is not smooth, then 
the singular locus of $V/G$ is of codim $\geq 2$, hence $\pi(W)$ is of codimension 2.

However, the singular set of $V/G$ is contained in $p(F)$ with $p: V \rightarrow V/G$ and hence $codim(F) \geq 2$.  It thus suffices to show that $Fix(g)$
is not of codimension one for any $g \in G$. But this can be excluded as any element $g \in G$ is symplectic.
\end{proof}

\begin{Rque}
The above proof shows that the same assertion holds if $V$ is only $\qit$-factorial with singularities in codimension $\geq 3$.
\end{Rque}

As in dimension 2, quotient singularities are A-D-E singularities which admit crepant resolutions, the above theorem shows:
\begin{Cor}
Let $V$ be a smooth irreducible symplectic variety and $G$ a finite group of symplectic automorphisms of $V$. Suppose that $V/G$ has only isolated
singularities. Then $V/G$ admits a symplectic resolution if and only if $dim(V) =2$.
\end{Cor}
\begin{Rque}
In \cite{SB}, N. Shepherd-Barron has proved that under some extra conditions, any symplectic resolution of an isolated symplectic singularity is isomorphic
to the collapsing of the zero-section in the cotangent bundle $T^*\pit^n \rightarrow \overline{\0}_{min}$, where $\overline{\0}_{min}$
is the minimal nilpotent orbit in $\mathfrak{sl}(n+1,\cit).$
\end{Rque}

\section{Case of $T^*X/G$}
In this section, we will restrict to the following case.
Let $X$ be an n-dimensional (quasi-projective) smooth irreducible complex variety and $G$ a
 finite group of automorphisms of $X$. Let $ T^*X$ be the total space of 
cotangent bundle of $X$ and $\omega$ the canonical symplectic form on $T^*X$.
 The group $G$ acts naturally on $T^*X$ preserving the symplectic form $\omega$,
thus $T^*X/G$ is a variety with symplectic singularities. Our purpose here is to give some necessary conditions on $G$ for the existence of a symplectic 
resolution for $T^*X/G$.
\begin{Thm}
If $T^*X/G$ admits a symplectic resolution, then $X/G$ is smooth.
\end{Thm} 

This theorem is essentially proved by Kaledin in \cite{Ka1} (and some of his other results extend also to our case), where he considered the
 special case of  $X = \cit^n$. Here we just give an outline of the proof. 

Consider the natural $\cit^*$-action on the vector bundle $T^*X = T^*X$, which commutes with the $G$-action, so we get a $\cit^*$-action on $T^*X/G$.
For this $\cit^*$-action, we have $\lambda^* \omega = \lambda \omega$. The fixed points of this action are identified with the subvariety 
 $X/G \subset T^*X/G$.
 Let $\pi: Z \rightarrow T^*X/G$ be a symplectic resolution and $\Omega$ the
corresponding symplectic form on $Z$. As shown in \cite{Ka1} (see also section 3.1 of \cite{Fu}), the $\cit^*$-action on $T^*X/G$ lifts to 
 $Z$ in  such a  way that  $\pi$ is $\cit^*$-equivariant. For this $\cit^*$-action, we have $\lambda^* \Omega = \lambda \Omega$.
The key lemma is 
\begin{Lem}
For every $x \in X/G \subset T^*X/G$, there exists at most  finitely many points in $\pi^{-1}(x)$ which are fixed by the $\cit^*$-action on $Z$.
\end{Lem}  

The proof is based on the equation $\lambda^* \Omega = \lambda \Omega$ and the semismallness of the map $\pi$. For details, see prop. 6.3 of \cite{Ka1}.
Since a generic point on $X/G$ is smooth, the map $\pi: \pi^{-1}(X/G) \rightarrow X/G$ is generically one-to-one and surjective, thus there exists a 
connected component $Y$ of fixed points $Z^{\cit^*}$ such that $\pi: Y \rightarrow X/G$ is dominant and generically one-to-one. Now by the above lemma,
 this map is also finite. Now that $X/G$ is normal implies that $\pi: Y \rightarrow X/G$ is in fact an isomorphism.
Since $Z$ is smooth, $Z^{\cit^*}$ is a union of smooth components, so $Y$ is smooth, thus $X/G$ is smooth. 

\begin{Thm}
If $T^*X/G$ admits a symplectic resolution $\pi: Z \rightarrow T^*X/G$, then $Z$ contains an open set $U$ which is isomorphic to $T^*(X/G)$.
\end{Thm}
\begin{proof}
By the proof of above theorem, there exists a connected component $Y$ of $Z^{\cit^*}$ such that 
$\pi: Y \rightarrow X/G$ is an isomorphism. In particular we have $dim(Y) = n$.

For any fixed point $y \in Y$, the action of $\cit^*$ on $Z$ induces a weight decomposition 
$$T_yZ = \oplus_{p \in \zit} T_y^p Z, $$ where $T_y^p Z = \{ v \in T_y Z |  \lambda_* v = \lambda^p v \}$, and $T_yY$ is identified to $T_y^0Z$.
The equation $\lambda^* \Omega = \lambda \Omega$ gives a duality between $T_y^p(Z)$
 and $T_y^{1-p}(Z)$. In particular,
$dim(T_y^1Z) = dim(T_y^0Z) = dim Y = n$, so $T_y^pZ =0$ for all $p \neq 0,1$, which gives a decomposition $T_yZ = T_yY \oplus T_y^1Z$ and that $Y$ is 
Lagrangian with respect to $\Omega$.

Let $U$ be the attraction subvariety of $Y$, i.e. $U = \{ z \in Z | \lim_{\lambda \rightarrow 0} \lambda \cdot z \in Y \}$, and $p: U \rightarrow Y$
the attraction map. By the work of A.Bialynicki-Birula \cite{BB}, the decomposition $T_yZ = T_yY \oplus T_y^1Z$ gives that $U$ a vector 
bundle of rank $n$ over $Y$, so $U$ is identified with the total space of the normal bundle $N$ of $Y$ in $Z$. 
Now we establish an isomorphism between $N$ and $T^*Y$ as follows. Denote by $\Omega_{can}$ the canonical symplectic structure on $T^*Y$.
Take a point $y \in Y$, and a vector $v \in N_y$.
Since $Y$ is Lagrangian in the both symplectic spaces, there exists a unique vector $w \in T_y^*Y$ such that
 $\Omega_y(v,u) = \Omega_{can,y}(w,u)$ for any $u \in T_yY$.
We define the map $i: N \rightarrow T^*Y$ to be $i(v)=w$. It is clear that $i$ is a $\cit^*$-equivariant isomorphism.
\end{proof}

Here are some immediate corollaries.
\begin{Cor}
If $X/G$ is simply connected, then $Z$ is also simply connected.
\end{Cor}
\begin{proof}
If $X/G$ is simply connected, so is $T^*(X/G)$. Notice that the open map $i: U \simeq T^*(X/G)$ induces a surjective map 
$\pi_1(U) \rightarrow  \pi_1(Z)$, so $\pi_1(Z) = 0$.
\end{proof}

Now we will consider the case of $X = \cit^n$ and $G$ a finite group in $GL(n, \cit)$. 
\begin{Cor}
Let $\pi: Z \rightarrow \cit^n\oplus \cit^n/G$ be a symplectic resolution, then $Z$ is simply connected and rational.
\end{Cor}
\begin{proof}
By our theorem, $\cit^n/G$ is smooth. A classical theorem (\S 5 Chap. V of \cite{Bo}) implies then $\cit^n/G$ is isomorphic to $\cit^n$, thus
$T^*(\cit^n/G)$ is simply connected and rational, which gives the corollary.
\end{proof}

That $Z$ is simply connected has been proved by M.Verbitsky in \cite{Ve}.
\begin{Cor}
Let $\pi: Z \rightarrow \cit^n\oplus \cit^n/G$ be a symplectic resolution and $F = \pi^{-1}(0)$ the exceptional set over 0.
 If $F$ is either $n$-dimensional irreducible or smooth, then $F$ is rational.
\end{Cor}
\begin{proof}
Denote by $y_0 \in Y$ the point mapped to $0 \in V\oplus V^*/G$. In both cases, $F$ is irreducible.
Since $\pi$ is $\cit^*$-equivariant, $F \cap U $ is contained in $T_{y_0}^*Y$.
If $F$ is of dimension $n$, then  $F \cap U$ is $n$-dimensional and closed in $U$, so $\cit^n \simeq T_{y_0}^*Y  \subset F$, i.e. $F$ is rational.
If $F$ is smooth, then $F \cap U$ is also smooth. Note that $F \cap U$ is a smooth closed cone in $\cit^n$, so $F \cap U$ is a linear space $\cit^k$,
with $k = dim(F)$, which gives that $F$ is rational.
\end{proof}

The sub-variety $F$ is the essential difficulty to understand $Z$. It looks like that every irreducible component of $F$ is rational, but we could not
prove this in general. For the 4-dimensional case, this has been proved by J.Wierzba \cite{Wi}. 

\section{McKay correspondence}

Let $X$ be an algebraic complex variety. The cohomology groups with compact supports $H_c^*(X, \cit)$ admit a canonical mixed Hodge structure.
We will denote by $h^{p,q}(H^k_c(X,\cit))$ the dimension of the $(p,q)$-Hodge component of the $k$-th cohomology.
\begin{Def}
We define $e^{p,q} : = \sum_{k \geq 0} (-1)^k h^{p,q}(H^k_c(X,\cit))$, and the polynomial $E(X;u,v):= \sum_{p,q} e^{p,q}(X) u^pv^q$ is called the 
{\em E-polynomial} of $X$.
\end{Def}

Here are some basic properties of $E$-polynomials, for proofs see \cite{BD}. \\
$\bullet$  If the Hodge structure is {\em pure}, i.e.  $h^{p,q}(H^k_c(X,\cit))= 0$ when  $p+q \neq k$, then $e^{p,q} = (-1)^{p+q} h^{p,q}(X)$. \\
$\bullet$  The Euler characteristic of $X$ is equal to $E(X; -1, -1)$. \\
$\bullet$  If $X$ is a disjoint union of locally closed subvarieties $X_i$, then $E(X;u,v) = \sum_i E(X_i; u,v)$. \\
$\bullet$  If $X \rightarrow Y$ is a locally trivial fibration with fiber $F$, then $E(X;u,v) = E(Y;u,v) \cdot E(F;u,v)$.  

Consider a smooth irreducible complex variety $V$. Let $G$ be a finite group in $Aut(V)$ and $g \in G$ an arbitrary element. 
To simplify the exposition, we suppose that $G$ preserves a volume form on $V$.
Take a connected component $W$ of $V^g$,
fixed points by $g$.
 For a point $x \in W$, the differential  $dg$ is an automorphism of $T_x V$.
 Denote by $e^{2\pi\sqrt{-1} a_i}$ the eigenvalues of 
$dg$ with $a_i \in \qit \cap [0,1[$.
The {\em weight $wt(g,W)$} is defined to be $\sum_i a_i$.
Since $G$ preserves a volume form, $dg$ has determinant 1, this gives that the weight is an integer, thus it  is indepent of the choice of the point $x$
in the connected component $W$.
 Moreover, if $h \in G$ commutes
with $g$, and $W' = h \cdot W$ another component of $X^g$, then $wt(g, W') = wt(g, W)$. Denote by $C(g)$ the centralizer of $g$ in $G$ and by
$C(g,W)$ the subgroup of all elements in $C(g)$ which leaves the component $W$ invariant. 
\begin{Def}
The {\em orbifold E-function} of a $G$-manifold $V$ is defined to be
 $$E_{orb}(V;G;u,v) := \sum_{\{g\}} \sum_{\{W\}} (uv)^{wt(g,W)} E(W/C(g,W); u,v), $$
where $\{g\}$ runs over all conjugacy classes in $G$, and $\{W\}$ runs over the set of representatives of all $C(g)$-orbits in the set of connected
components of $V^g$.
\end{Def}

In the paper \cite{Ba}, V.Batyrev has proved that the orbitfold E-function is equal to the string E-function. Applying this to our case, we get
\begin{Thm} \label{thm2}
Let $X$ be an $n$-dimensional smooth projective variety and $G$ a finite group in $Aut(X)$.
Suppose there exists a symplectic resolution $\pi: Z \rightarrow T^*X/G$. Then: (1) the Hodge structure on $Z$ is pure;
(2) the Hodge numbers of $Z$ can be calculated by the following:
 $$\sum_{p,q} (-1)^{p+q} h^{p,q}(Z) u^p v^q = (uv)^n \sum_{\{g\}} \sum_{\{W\}}  E(W/C(g,W); u,v) , $$
where $\{g\}$ runs over all conjugacy classes in $G$ and 
$\{W\}$ runs over the set of representatives of all $C(g)$-orbits in the set of connected components of $X^g$.
\end{Thm}
\begin{proof}
Consider the natural  $\cit^*$-action on $T^*X/G$, which lifts to $Z$. We identify $X/G$ with the fixed points by the $\cit$-action in 
 $T^*X/G$.
 Since $X/G$ is projective and the map $\pi: \pi^{-1}(X/G) \rightarrow X/G$
is projective,  the variety $\pi^{-1}(X/G)$ is projective. Notice that the fixed points $Z^{\cit^*}$ by the $\cit^*$-action on $Z$ is contained
in $\pi^{-1}(X/G)$ and $Z$ is smooth, 
so  $Z^{\cit^*}$ is a union of smooth projective varieties. In particular, the Hodge structure on $Z^{\cit^*}$ is pure.
Now the same argument of the proof of theorem 8.4 in \cite{Ba} gives the affirmation (1).

For the second affirmation, theorem 7.5 of {\em loc.cit.} gives that $E(Z;u,v) = E_{orb}(T^*X; G; u,v)$.
Since the Hodge structure on $Z$ is pure, $E(Z;u,v) =\sum_{p,q} (-1)^{p+q} h^{p,q}(Z) u^p v^q $, thus we need to calculate $E_{orb}(T^*X; G; u,v)$.
The key point here is that every fixed component $F$ in $T^*X$ is of the form $T^*W$, for some fixed component $W$ in $X$ and $C(g,F) = C(g,W)$. 
The map $F/C(g,F) \rightarrow W/C(g, W)$ is a fibration with fiber $\cit^{dim(W)}$, so 
$$E(F/C(g,F);u,v) =  (uv)^{dim(W)}  E(W/C(g,W);u,v).$$

Now what we need is to calculate the weight of $g$ at $W$. Since $g$ acts on $T^*X$ as a symplectic automorphism, lemma 2.6 of \cite{Ka2} gives that 
$wt(g, W) = \frac{1}{2} codim(T^*W) =  codim(W)$. Now our theorem follow immediately. 
\end{proof}

A particular simple situation is the following:
\begin{Thm}
Let $X= \pit^n$ be the projective space and $G$ a finite group in $SL(n+1,\cit)$. Suppose that $\pi: Z \rightarrow T^*X/G $ is a symplectic
resolution. Then the Hodge numbers of $Z$ can be calculated by the following 
$$\sum_{p,q} (-1)^{p+q} h^{p,q}(Z) u^p v^q = (uv)^n \sum_{\{g\}} \sum_i \cfrac{(uv)^{k_i(g)} - 1}{uv - 1} ,  $$ where
$k_i(g)$ are multiplicities of distinct eigenvalues of $g \in SL(n+1,\cit)$ and $\{g\}$ runs over conjugacy classes
 in $G$, which is considered as a subgroup in $Aut(\pit^n)$.
\end{Thm}
\begin{proof}
Let $p: \cit^{n+1} - \{0\} \rightarrow \pit^n$ be the natural projection. 
For any $g \in SL(n+1,\cit)$,  the point $p(x) \in \pit^n$ is fixed by $g$ if and only if $g(x) = \lambda x$ for some $\lambda$.
In particular, $x$ should be an eigenvector of $g$.
 If we denote by $L_i$ the eigenspaces
of $g$,  then the  fixed points $X^g$ are  $X_i^g = p(L_i)$. Notice that if $h \in G$ commutes with $g$, then 
$g(h(x)) = \lambda h(x)$, i.e. $h(x)$ and $x$ are in the same eigenspace of $g$, which shows that the $C(g)$-action on the set of connected
 components of $X^g$ is trivial.
 So the index set $\{W\}$ in the sum of the orbit  E-function is the same as $\{ X_i^g\}$.

Consider an eigenspace $L_i$, whose dimension is denoted by  $k_i(g)$.
  Then there is a fibration $L_i - \{0\}/C(g) \rightarrow X_i^g/C(g)$ with
fiber $\cit^*$. Now use properties of E-polynomials, we have $E(L_i - \{0\}/C(g); u,v) = E(L_i/C(g); u,v) - E(\{0\}; u,v) = (uv)^{k_i(g)} - 1 $
and $E(\cit^*; u,v) = uv - 1$. This gives that $E(X_i^g/C(g); u, v) = \cfrac{(uv)^{k_i(g)} - 1}{uv - 1}$. Now the theorem follows directly 
from theorem \ref{thm2}.
\end{proof}

\begin{Cor}
The Euler number of $Z$ is given by $e(Z) = (n+1) c(G) $, where $c(G)$ is the number of conjugacy classes of elements in $G$, which is considered as a
subgroup of $Aut(\pit^n)$.
\end{Cor}
\begin{Cor}
 $H^{2j+1}(Z,\cit) = 0$  for all $j$ and $H^{2i}(Z,\cit) = 0$ for $i \leq n-1$. Furthermore $H^{2i}(Z,\cit)$  has the Hodge type $(i,i)$ for all $i$.
\end{Cor}

\section{One application}
 Let $\Sigma$ be a 
smooth projective curve. There is a natural action of the permutation group $S_n$ on the smooth projective   variety $\Sigma^n$.
This action gives an action of $S_n$ on $T^*(\Sigma^n) \simeq (T^*\Sigma)^n$. As easily seen, the action of $S_n$ on $(T^*\Sigma)^n$
is exactly the permutation action, so $T^*(\Sigma^n)/S_n \simeq S^n(T^*\Sigma)$. 
The well-known  resolution $Hilb^n(T^*\Sigma) \rightarrow S^n(T^*\Sigma)$ gives a symplectic resolution for the quotient $T^*(\Sigma^n)/S_n$.
Now we will use the McKay correspondence to deduce the cohomology of $Hilb^n(T^*\Sigma)$.

Notice that the conjugacy classes of the permutation group $S_n$ correspond to partitions of $n$. Let $\nu = (1^{a_1} 2^{a_2} \cdots n^{a_n})$
be a partition of $n$, i.e. $a_i \geq 0$ and $\sum_j ja_j = n$. Here $a_i$ is the number of length $i$-cycles in the permutation. 
 Then the fixed points subvariety $Fix(\nu)$ in $\Sigma^n$ by $\nu$ is isomorphic to 
$\prod_j \Sigma^{a_j}$. This gives that $Fix(\nu)/C(\nu) \simeq S^\nu \Sigma : = \prod_j S^{a_j} \Sigma$. By our theorem \ref{thm2}, we have 
$$\sum_k h^{2k}(Hilb^n(T^*\Sigma)) t^{2k} = t^{2n} \sum_{\nu } \sum_j h^{2j}(S^\nu \Sigma) t^{2j}.$$
By the Poicar\'e duality and change $t$ to $1/t$,  this gives 
\begin{Lem}
$$ P_t(Hilb^n(T^*\Sigma)) = \sum_\nu t^{2n - 2d(\nu)} P_t(S^\nu \Sigma),$$
where $d(\nu) = \sum_j a_j$ is the complex dimension of $S^\nu \Sigma$.
\end{Lem}
Now using Macdonald's formula (see \cite{Ma}) for cohomology of symmetric products of algebraic curves, 
we can deduce G\"ottsche's formula in the case of $T^*\Sigma$:
\begin{Prop}
$$\sum_{n=0}^\infty P_t(Hilb^n(T^*\Sigma)) q^n = \prod_{d=1}^\infty \cfrac{(1+t^{2d-1}q^d)^{b_1(\Sigma)}}{(1- t^{2d-2}q^d)^{b_0(\Sigma)}
(1- t^{2d}q^d)^{b_2(\Sigma)}}. $$
\end{Prop}
\begin{Rque}
In chap.7 \cite{Nak}, H.Nakajima has deduced this formula via an elementary argument.
He firstly proved that the fixed points subvariety of the $\cit^*$-action in $Hilb^n(T^*\Sigma)$ is the disjoint union $\coprod_\nu S^\nu \Sigma$, 
then using Morse theory, he proved lemma 7.1.
\end{Rque}

\quad \\
Labortoire J.A.Dieudonn\'e, Parc Valrose \\ 06108 Nice cedex 02, FRANCE \\
baohua.fu@polytechnique.org

\begin{thebibliography}{10}
\bibitem[Ba]{Ba}
V. Batyrev, \emph{Non-Archimedean integrals and stringy Euler numbers of log-terminal pairs}, J. Eur. Math. Soc. {\bf 1}, no. 1, 5-33(1999). 
\bibitem[BD]{BD}
V. Batyrev and D. Dais, \emph{Strong McKay correspondence, string-theoretic Hodge numbers and mirror symmetry}, Topology {\bf 35}, No.4, 901-929 (1996).
\bibitem[Be]{Be}
A. Beauville, \emph{Symplectic singularities}, Invent. Math. {\bf 139}, 541-549(2000).
\bibitem[BB]{BB}
A. Bialynicki-Birula, \emph{Some theorems on actions of algebraic groups},  Ann. of Math., II. Ser. {\bf 98}, 480-497 (1973).
\bibitem[Bo]{Bo}
N. Bourbaki, \emph{Groupes et alg\`ebres de Lie}, Chap.IV--VI, Hermann(Paris) 1968.
\bibitem[Fu]{Fu}
B. Fu, \emph{Symplectic resolutions for nilpotent orbits}, preprint math.AG/0205048, to appear in Invent. Math.
\bibitem[Ka1]{Ka1}
D. Kaledin, \emph{Dynkin diagrams and crepant resolutions of quotient singularities}, preprint math.AG/9903157.
\bibitem[Ka2]{Ka2}
D. Kaledin, \emph{McKay correspondence for symplectic quotient singularities}, Invent. Math. {\bf 148} 1, 150-175(2002).
\bibitem[Ka3]{Ka3}
D. Kaledin, \emph{Symplectic resolutions: deformations and birational maps}, preprint math.AG/0012008.
\bibitem[Ma]{Ma}
I. Macdonald, \emph{The Poicar\'e polynomial of a symmetric product}, Proc. Camb. Phil. Soc. {\bf 58}, 563-568 (1962).
\bibitem[Nak]{Nak}
H. Nakajima, \emph{Lectures on Hilbert schemes of points on surfaces}, University Lecture Series {\bf 18}, Providence AMS, 1999.
\bibitem[Nam]{Nam}
Y. Namikawa, \emph{Deformation theory of singular symplectic n-folds}, Math. Ann. {\bf 319}, no. 3, 597-623 (2001).
\bibitem[Re]{Re}
M. Reid, \emph{La correspondance de McKay}, S\'eminaire Bourbaki n.867, preprint math.AG/9911165.
\bibitem[SB]{SB}
N. Shepherd-Barron, \emph{Long extremal rays and symplectic resolutions}, preprint. 
\bibitem[Ve]{Ve}
M. Verbitsky,  \emph{Holomorphic symplectic geometry and orbifold singularities}, Asian J. Math. {\bf 4}, no. 3, 553-563 (2000).
\bibitem[Wi]{Wi}
J. Wierzba, \emph{Symplectic Singularities}, Ph.D. thesis, Trinity college, Cambridge University, September 2000.
\end{thebibliography}
\end{document}